\pdfoutput=1 

\documentclass[a4paper]{amsart}

\usepackage[british]{babel}

\usepackage{microtype}  
\usepackage{mathrsfs}   
\usepackage{upgreek}    
\usepackage{hyperref}   
\usepackage{graphicx}   

\newtheorem{theorem}{Theorem}           
\newtheorem*{proposition}{Proposition}  
\theoremstyle{definition}   
\newtheorem*{definition}{Definition}    
\newtheorem{example}{Example}           

\DeclareMathOperator{\with}{:}          
\DeclareMathOperator{\restrict}{\llcorner}   
\DeclareMathOperator{\Clos}{Closure}    
\DeclareMathOperator{\Tan}{Tan}         
\DeclareMathOperator{\spt}{spt}         
\DeclareMathOperator{\trace}{trace}     
\DeclareMathOperator{\image}{image}     
\DeclareMathOperator{\Lip}{Lip}         
\DeclareMathOperator{\grad}{grad}       
\DeclareMathOperator{\Nor}{Nor}         
\DeclareMathOperator{\Der}{D}           

\newcommand{\ud}{\ensuremath{\,\mathrm{d}}}

\newcommand{\version}[2]{#2}         

\newenvironment{pullquote}
	{\begin{figure}[ht]
		\rule{\textwidth}{1pt}
		\begin{center}
			\noindent \hspace{\stretch{1}}}
	{		\hspace{\stretch{1}}
		\end{center}
		\rule[3pt]{\textwidth}{1pt}
	\end{figure}}

\title{The concept of varifold}
\author{Ulrich Menne}

\version{\thanks{Ulrich Menne is a post-doctoral researcher at the University
of Zurich.  His e-mail address is \nolinkurl{Ulrich.Menne@uzh.ch}.}}
{\thanks{The author is grateful to Prof.~Frank Morgan for his detailed
comments on this note.}}

\version{}{\date{\today}}

\begin{document}



\version{}{\address{Department of Mathematics, Faculty of Science, University
of Zurich, Win\-ter\-thu\-rer\-stras\-se 190, \textsc{8057 Zurich,
Switzerland}} \email{Ulrich.Menne@uzh.ch}}

\begin{abstract}
	We survey -- by means of 20~examples -- the concept of varifold, as
	generalised submanifold, with emphasis on regularity of integral
	varifolds with mean curvature, while keeping prerequisites to a
	minimum.  Integral varifolds are the natural language for studying the
	variational theory of the area integrand if one considers, for
	instance, existence or regularity of stationary (or, stable) surfaces
	of dimension at least three, or the limiting behaviour of sequences of
	smooth submanifolds under area and mean curvature bounds.
\end{abstract}


\maketitle


\section{Introduction}

\subsection*{Motivation}

Apart from generalisation, there are two main reasons to include non-smooth
surfaces in geometric variational problems in Euclidean space or Riemannian
manifolds: firstly, the separation of existence proofs from regularity
considerations and, secondly, the modelling of non-smooth physical objects.
Following the first principle, varifolds were introduced \version{by
F.~Almgren }{}in~1965, to prove, for every intermediate dimension, the
existence of a generalised \emph{minimal surface} (i.e., a surface with
vanishing first variation of area) in a given compact smooth Riemannian
manifold\version{}{ (see~\cite{Almgren:Vari})}.  Then, in~1972, an important
partial regularity result for such varifolds was established \version{by
W.~Allard in his foundational paper \emph{On the first variation of a
varifold}}{in}~\cite{MR0307015}.  These pioneering works still have a
strong influence in geometric analysis as well as related fields.
Particularly prominent examples are the proof of the Willmore conjecture
\version{by F.~Marques and A.~Neves}{(see~\cite{marques_neves_willmore})} and
the development of a regularity theory for stable generalised minimal surfaces
in one codimension \version{by N.~Wickramasekera}{(see~\cite{MR3171756})},
both published in~2014.  Following the second principle, varifolds were
employed\version{ by Almgren's PhD student K.~Brakke}{} in~1978 to create a
mathematical model of the motion of grain boundaries in an annealing pure
metal\version{}{ (see~\cite{MR485012})}.  This was the starting point of the
rapid development of mean curvature flow, even for smooth surfaces.

\subsection*{What is a varifold?}

Constructing non-parametric models of non-smooth surfaces usually requires
entering the realm of geometric measure theory.  In the simplest case, we
associate, to each smooth surface~$M$, the measure \emph{over the ambient
space} whose value at a set~$A$ equals the surface measure of the
intersection~$A \cap M$.  For varifolds, it is in fact expedient to
similarly record information on the tangent planes of the surfaces.  This
yields weak continuity of area, basic compactness theorems, and the
possibility to retain tangential information in the limit.

\subsection*{Main topics covered}

We will introduce the notational infrastructure together with a variety of
examples.  This will allow us to formulate the compactness theorem for
integral varifolds, Theorem~\ref{thm:integral_compactness}, and two key
regularity theorems for integral varifolds with mean curvature, Theorems
\ref{thm:allards_regularity}~and~\ref{thm:second_order_rectifiability}.
\version{The accompanying drawings generally stress certain aspects deemed
important at the expense of accuracy.  A version of this article with
additional references is available from
\url{https://arxiv.org/abs/1705.05253}.}{The present text is a version with
additional references but without figures of a note compiled for the
\emph{Notices of the American Mathematical Society}.}

\subsection*{Notation}

\emph{Suppose throughout this note, that $m$ and $n$ are integers and $1 \leq
m \leq n$}.  We make use of the $m$-di\-men\-sion\-al Hausdorff
measure~$\mathscr H^m$ over~$\mathbf R^n$.  This measure is one of several
natural outer measures over $\mathbf R^n$ which associate the usual values
with subsets of $m$-di\-men\-sion\-al continuously differentiable submanifolds
of~$\mathbf R^n$.


\section{First order quantities}

\version{\begin{pullquote}
	the attempt to comprise all objects that should be considered
	$m$-di\-men\-sion\-al surfaces
\end{pullquote}}{}

Amongst the reasons of developing the notion of \emph{$m$-di\-men\-sion\-al
varifold} in $\mathbf R^n$ was the attempt to comprise all objects that
should be considered $m$-di\-men\-sion\-al surfaces of locally finite area
in~$\mathbf R^n$.  Its definition is tailored for compactness.

\begin{definition} [Varifold and weight\version{}{,
see~\protect{\cite[3.1]{MR0307015}}}] \label{def:varifold}

	If $\mathbf G (m,n)$ is the space of unoriented $m$-di\-men\-sion\-al
	subspaces of~$\mathbf R^n$, endowed with its natural topology, then,
	by an \emph{$m$-di\-men\-sion\-al varifold~$V$} in $\mathbf R^n$, we
	mean a Radon measure $V$ over $\mathbf R^n \times \mathbf G (n,m)$,
	and we denote by~$\|V\|$, the \emph{weight} of~$V$, that is, its
	canonical projection into $\mathbf R^n$.
\end{definition}

The theory of Radon measures\version{}{ (see~\cite[2.23]{MR3528825})} yields a
metrisable topology on the space of $m$-di\-men\-sion\-al varifolds in
$\mathbf R^n$ such that $V_i \to V$ as $i \to \infty$ if and only if
\begin{equation*}
	{\textstyle \int k \ud V_i \to \int k \ud V} \quad \text{as $i \to
	\infty$}
\end{equation*}
whenever $k : \mathbf R^n \times \mathbf G (m,n) \to \mathbf R$ is a
continuous function with compact support.

\begin{proposition} [Basic compactness\version{}{,
see~\protect{\cite[2.6\,(2a)]{MR0307015}}}] \label{thm:compactness}

	Whenever $\kappa$ is a real-valued function on the bounded open
	subsets of~$\mathbf R^n$, the set of $m$-di\-men\-sion\-al
	varifolds~$V$ in~$\mathbf R^n$ satisfying $\| V \| ( Z ) \leq \kappa
	(Z)$ for each bounded open subset~$Z$ of~$\mathbf R^n$, is compact.
\end{proposition}

In general, neither does the weight~$\| V \|$ of~$V$ need to be
$m$~di\-men\-sion\-al, nor does the Grassmannian information of~$V$ need to be
related to the geometry of~$\| V \|$:

\begin{example} [Region with associated plane] \label{example:region}
	For each plane $T \in \mathbf G (n,m)$ and each open subset~$Z$
	of~$\mathbf R^n$, the product of, the Lebesgue measure restricted to
	$Z$, with the Dirac measure concentrated at~$T$, is an
	$m$-di\-men\-sion\-al varifold in~$\mathbf R^n$.
\end{example}

\begin{example} [Point with associated plane] \label{example:point}
	Whenever $a \in \mathbf R^n$ and $T \in \mathbf G (n,m)$, the Dirac
	measure concentrated at~$(a,T)$ is an $m$-di\-men\-sion\-al varifold
	in $\mathbf R^n$.
\end{example}

This generality is expedient to describe higher (or, lower) dimensional
approximations of more regular $m$-di\-men\-sion\-al surfaces, so as to
include information on the $m$-di\-men\-sion\-al tangent planes of the limit.
Such situations may be realised by elaborating on
Examples~\ref{example:region} and~\ref{example:point}.  However, we will
mainly consider \emph{integral} varifolds; that is, varifolds in the spirit of
the next basic example, which avoid the peculiarities of the preceding two
examples.

\begin{example} [Part of a submanifold\version{}{,
see~\protect{\cite[3.5]{MR0307015}}}] \label{example:part}

	If $B$ is a Borel subset of a closed $m$-di\-men\-sion\-al
	continuously differentiable submanifold $M$ of $\mathbf R^n$, then, by
	Riesz's representation theorem, the associated varifold~$\mathbf v_m (
	B )$ in~$\mathbf R^n$ may be defined by
	\begin{equation*}
		{\textstyle \int k \ud \mathbf v_m ( B ) = \int_B k ( x, \Tan
		(M,x) ) \ud \mathscr H^m \, x}
	\end{equation*}
	whenever $k : \mathbf R^n \times \mathbf G (m,n) \to \mathbf R$ is a
	continuous function with compact support, where $\Tan (M,x)$ is the
	tangent space of $M$ at $x$.  Then, the weight~$\| \mathbf v_m (B) \|$
	over~$\mathbf R^n$ equals the restriction~$\mathscr H^m \restrict B$
	of $m$~di\-men\-sion\-al Hausdorff measure to~$B$.
\end{example}

\begin{definition} [Integral varifold\version{}{,
see~\protect{\cite[3.5]{MR0307015}}}] \label{def:integral_varifold}

	An $m$-di\-men\-sion\-al varifold~$V$ in~$\mathbf R^n$ is
	\emph{integral} if and only if there exists a sequence of Borel
	subsets~$B_i$ of closed $m$-di\-men\-sion\-al continuously
	differentiable submanifolds $M_i$ of~$\mathbf R^n$ with $V =
	\sum_{i=1}^\infty \mathbf v_m (B_i)$.
\end{definition}

An integral varifold~$V$ is determined by its weight~$\| V \|$; in fact,
\version{}{see \cite[3.5\,(1)]{MR0307015}, }associated to~$\| V \|$, there are
a multiplicity function~$\Theta$ and a tangent plane function~$\tau$ with
values in the nonnegative integers and $\mathbf G(n,m)$, respectively, such
that \begin{equation*} {\textstyle V(k) = \int k (x,\tau(x)) \Theta (x) \ud
\mathscr H^m \,x} \end{equation*} whenever $k : \mathbf R^n \times \mathbf G
(m,n) \to \mathbf R$ is a continuous function with compact support.

\begin{example} [Sum of two submanifolds\version{}{,
see~\protect{\cite[3.5\,(2)]{MR0307015}}}] \label{example:sum}

	Suppose $M_i$ are closed $m$-di\-men\-sion\-al continuously
	differentiable submanifolds~$M_i$ of~$\mathbf R^n$, for $i \in \{ 1,2
	\}$, and $V = \mathbf v_m (M_1) + \mathbf v_m (M_2)$.  Then, for
	$\mathscr H^m$ almost all $z \in M_1 \cap M_2$, we have
	\begin{equation*}
		\Theta (z) = 2 \quad \text{and} \quad \tau (z) = \Tan (M_1,z)
		= \Tan (M_2,z).
	\end{equation*}
\end{example}

The theory of varifolds is particularly successful in the variational study of
the area integrand.  To describe some aspects of this, diffeomorphic
deformations are needed in order to define and compute its first variation.

\begin{definition} [Area of a varifold]
	If $B$ is a Borel subset of $\mathbf R^n$, and $V$ is an
	$m$-di\-men\-sion\-al varifold $V$ in $\mathbf R^n$, then \emph{the
	area of~$V$ in $B$} equals $\| V \| (B)$.
\end{definition}

\begin{example} [Area via multiplicity\version{}{,
see~\protect{\cite[3.5\,(1)]{MR0307015}}}]
\label{example:area_via_multiplicity}

	If $V$ is integral, then $\| V \| (B) = \int_B \Theta \ud \mathscr
	H^m$.
\end{example}

\begin{definition} [Push forward\version{}{,
see~\protect{\cite[3.2]{MR0307015}}}] \label{def:push_forward}

	Suppose $V$~is an $m$-di\-men\-sion\-al varifold in~$\mathbf R^n$, and
	$\phi$~is a proper continuously differentiable diffeomorphism
	of~$\mathbf R^n$ onto $\mathbf R^n$.  Then, the \emph{push forward~of
	$V$ by~$\phi$} is the $m$-di\-men\-sion\-al varifold~$\phi_\# V$
	defined by
	\begin{equation*}
		{\textstyle \int k \ud \phi_\# V = \int k ( \phi (z),
		\image ( \Der \phi (z) |S) ) j_m \phi (z,S) \ud V \, (z,S)}
	\end{equation*}
	whenever $k : \mathbf R^n \times \mathbf G (m,n) \to \mathbf R$ is a
	continuous function with compact support; here, $j_m \phi (z,S)$
	denotes the $m$-di\-men\-sion\-al Jacobian of the restriction to~$S$
	of the differential $\Der \phi (z) : \mathbf R^n \to \mathbf R^n$.
\end{definition}

\begin{example} [Push forward\version{}{,
see~\protect{\cite[3.5\,(3)]{MR0307015}}}] \label{example:push_forward}

	If $B$ and $M$ are as in Example~\ref{example:part}, then, by the
	transformation formula (or area formula), $\phi_\# ( \mathbf v_m (B))
	= \mathbf v_m ( \image (\phi |B) )$.
\end{example}

\begin{definition} [First variation\version{}{,
see~\protect{\cite[4.2]{MR0307015}}}] \label{def:first_variation}

	The \emph{first variation} of a varifold~$V$ in~$\mathbf R^n$ is
	defined by
	\begin{equation*}
		( \updelta V ) ( \theta ) = {\textstyle\int \trace ( \Der
		\theta (z) \circ S_\natural ) \ud V (z,S) }
	\end{equation*}
	whenever $\theta : \mathbf R^n \to \mathbf R^n$ is a smooth
	vector field with compact support, where $S_\natural$~denotes the
	canonical orthogonal projection of~$\mathbf R^n$ onto~$S$.
\end{definition}

\begin{example} [Meaning of the first variation\version{}{,
see~\protect{\cite[4.1]{MR0307015}}}] \label{example:diffeomorphisms}

	If $\phi_t : \mathbf R^n \to \mathbf R^n$ satisfy $\phi_t (z) =
	z+t\theta (z)$ whenever $t \in \mathbf R$, $|t| \Lip \theta < 1$, and
	$z \in \mathbf R^n$, then, one may verify, that, for such~$t$,
	$\phi_t$ is a proper smooth diffeomorphism of~$\mathbf R^n$ onto
	$\mathbf R^n$, and $\| (\phi_t )_\# V \| ( \spt \theta )$ is a smooth
	function of $t$, whose differential at $0$ equals $( \updelta V ) (
	\theta )$.
\end{example}

\begin{example} [First variation of a smooth submanifold I\version{}{,
see~\protect{\cite[4.4, 4.7]{MR0307015}}}]
\label{example:variation_submanifold}

	Suppose $M$ is a properly embedded, twice continuously differentiable
	$m$-di\-men\-sion\-al submanifold-with-boundary in~$\mathbf R^n$, and
	$V = \mathbf v_m (M)$.  Then, integrating by parts, one obtains
	\begin{equation*}
		{\textstyle ( \updelta V ) ( \theta ) = - \int_M \mathbf h
		(M,z) \bullet \theta (z) \ud \mathscr H^m \, z +
		\int_{\partial M} \nu (M,z) \bullet \theta (z) \ud \mathscr
		H^{m-1} \, z},
	\end{equation*}
	where $\mathbf h (M,z) \in \mathbf R^n$~is the mean curvature vector
	of~$M$ at~$z \in M$, $\nu (M,z)$~is the outer normal of~$M$ at~$z \in
	\partial M$, and $\bullet$~denotes the inner product in~$\mathbf R^n$.
	In fact, $\mathbf h(M,z)$ belongs to the normal space, $\Nor (M,z)$,
	of $M$ at~$z$.
\end{example}

\section{Second order quantities}

\version{\begin{pullquote}
	for varifolds of locally bounded first variation, one may define mean
	curvature.
\end{pullquote}}{}

In general, even integral varifolds -- which are significantly more regular
than general varifolds -- give meaning only to \emph{first} order properties
such as tangent planes, but do not possess any \emph{second} order
properties such as curvatures.  Yet, for varifolds of locally bounded first
variation, one may define mean curvature.

\begin{example} [A continuously differentiable submanifold without -- even
approximate -- curvatures, see\version{ R.~Kohn
}{~}\protect{\cite{MR0427559}}] \label{example:kohn} 

	There exists a nonempty closed continuously differentiable
	$m$-di\-men\-sion\-al submanifold~$M$ of~$\mathbf R^{m+1}$ such that,
	for any \emph{twice} continuously differentiable $m$-di\-men\-sion\-al
	submanifold~$N$ of~$\mathbf R^{m+1}$, there holds $\mathscr H^m ( M
	\cap N ) = 0$.
\end{example}

\begin{definition} [Locally bounded first varation\version{}{,
see~\protect{\cite[39.2]{MR756417}}}]
\label{def:locally_bounded_first_variation}

	A varifold~$V$ in~$\mathbf R^n$ is said to be \emph{of locally bounded
	first variation} if and only if $\| \updelta V \| (Z)$, defined as
	\begin{equation*}
		\sup \{ ( \updelta V ) (\theta) \with
		\text{$\theta$ a smooth vector field compactly supported in
		$Z$ with $| \theta | \leq 1$} \},
	\end{equation*}
	whenever $Z$ is an open subset of~$\mathbf R^n$, is finite on bounded
	sets.  Then, $\| \updelta V \|$~may be uniquely extended to a Radon
	measure over~$\mathbf R^n$, also denoted by~$\| \updelta V \|$.
\end{definition}

Clearly, for any open subset~$Z$ of~$\mathbf R^n$, $\| \updelta V \| (Z)$ is
lower semicontinuous in $V$.

\begin{example} [Functions of locally bounded variation\version{}{,
see~\protect{\cite[4.5]{MR0307015}}}]

	Suppose $m = n$.  Then, $m$-di\-men\-sion\-al varifolds in~$\mathbf
	R^m$ of locally bounded first variation are in natural correspondence
	to real valued functions of locally bounded first variation
	on~$\mathbf R^m$.
\end{example}

\begin{example} [First variation of a smooth submanifold II\version{}{,
see~\protect{\cite[4.4, 4.7]{MR0307015}}}]

	If $M$ and $V$ are as in Example~\ref{example:variation_submanifold},
	then $\| \updelta V \| = | \mathbf h (M,\cdot) | \mathscr H^m
	\restrict M + \mathscr H^{m-1} \restrict \partial M$.
\end{example}

\begin{theorem} [Compactness of integral varifolds, see\version{ Allard
}{~}\protect{\cite[6.4]{MR0307015}}] \label{thm:integral_compactness}

	If $\kappa$ is a real-valued function on the bounded open subsets
	of~$\mathbf R^n$, then the set of integral varifolds with $\| V \| ( Z
	) + \| \updelta V \| (Z) \leq \kappa (Z)$ for each bounded open
	subset~$Z$ of~$\mathbf R^n$, is compact.
\end{theorem}

The preceding theorem is fundamental both for varifolds and the consideration
of limits of submanifolds.  The summand $\| \updelta V \| (Z)$ therein may not
be omitted:

\begin{example} [Lattices of spheres\version{, see
Figure~\ref{fig:diffusing}}{}] \label{example:diffusing}

	If $V_i$ is the varifold associated to
	\begin{equation*}
		M_i = \mathbf R^{m+1} \cap \big \{ z \with \text{$|iz-a| =
		3^{-1} i^{-1-1/m}$ for some $a \in \mathbf Z^{m+1}$} \big \},
	\end{equation*}
	for every positive integer~$i$, then the non-zero limit of this
	sequence is the product the Lebesgue measure over~$\mathbf R^{m+1}$
	with an invariant Radon over $\mathbf G(m+1,m)$.
	\version{\begin{figure}
		\includegraphics[height=1.5in]{Lattice_of_spheres_big.jpg}
		\includegraphics[height=1.5in]{Lattice_of_spheres_small.jpg}
		\caption{One-di\-men\-sion\-al submanifolds (drawn in blue)
		of~$\mathbf R^2$ converging to a limit non-zero, non-integral
		varifold with all one-di\-men\-sion\-al tangent planes equally
		weighted at every point of~$\mathbf R^2$.}
		\label{fig:diffusing}
	\end{figure}}{}
\end{example}

\begin{definition} [Mean curvature\version{}{,
see~\protect{\cite[4.3]{MR0307015}}}] \label{def:mean_curvature}

	Suppose $V$ is a varifold in~$\mathbf R^n$ of locally bounded first
	variation, and $\sigma$ is the singular part of~$\| \updelta V \|$
	with respect to $\| V \|$.  Then, there exist a $\sigma$~almost
	unique, $\sigma$~measurable function $\boldsymbol{\eta} (V,\cdot)$
	with values in the unit sphere~$\mathbf S^{n-1}$ and a $\| V
	\|$~almost unique, $\| V \|$~measurable, locally $\| V\|$~summable,
	$\mathbf R^n$~valued function $\mathbf {h} (V,\cdot)$ -- called the
	\emph{mean curvature of~$V$} -- satisfying
	\begin{equation*}
		{\textstyle ( \updelta V ) ( \theta ) = - \int \mathbf h (V,z)
		\bullet \theta (z) \ud \| V \| \, z + \int \boldsymbol{\eta}
		(V,z) \bullet \theta (z) \ud \sigma \, z}
	\end{equation*}
	whenever $\theta : \mathbf R^n \to \mathbf R^n$ is a smooth
	vector field with compact support.
\end{definition}

\begin{example} [Varifold mean curvature of a smooth submanifold\version{}{,
see~\protect{\cite[4.4, 4.7]{MR0307015}}}] \label{example:varifold_mean_curv}

	If $M$, $V$, and~$\nu$ are as in
	Example~\ref{example:variation_submanifold}, then, $\sigma = \mathscr
	H^{m-1} \restrict \partial M$, $\boldsymbol{\eta} (V,z) = \nu(M,z)$
	for $\sigma$~almost all $z$, and $\mathbf h (V,z) = \mathbf h (M,z)$
	for $\| V \|$~almost all $z$.
\end{example}

Yet, in general, the mean curvature vector may have tangential parts, and the
behaviour of $\sigma$ and $\boldsymbol{\eta} (V,\cdot)$ may differ from that
of boundary and outer normal:

\begin{example} [Weighted plane\version{}{,
see~\protect{\cite[7.6]{MR3528825}}}]

	If $m=n$, $\Theta$ is a positive, continuously differentiable
	function on $\mathbf R^m$, and $V(k) = \int k (z,\mathbf R^m) \Theta
	(z) \ud \mathscr H^m \, z$ whenever $k : \mathbf R^m \times \mathbf G
	(m,m) \to \mathbf R$ is a continuous function with compact support,
	then
	\begin{equation*}
		\mathbf h (V,z) = \grad ( \log \circ \, \Theta ) (z) \quad
		\text{for $\| V\|$~almost all~$z$}.
	\end{equation*}
\end{example}

\begin{example} [Primitive of the Cantor function\version{}{,
see~\protect{\cite[12.3]{MR3528825}}}] \label{example:primitive}

	If $C$ is the Cantor set in $\mathbf R$, $f : \mathbf R \to \mathbf R$
	is the associated function (i.e., $f(x) = \mathscr H^d ( C \cap \{
	\chi \with \chi \leq x \} )$ for $x \in \mathbf R$, where $d = \log
	2/\log 3$), and $V$ is the varifold in~$\mathbf R^2 \simeq \mathbf R
	\times \mathbf R$ associated to the graph of a primitive function
	of~$f$, then, $V$ is an integral varifold of locally bounded first
	variation, $\mathbf h (V,z) = 0$ for $\| V \|$~almost all~$z$, and
	$\spt \sigma$ corresponds to~$C$ via the orthogonal projection onto
	the domain of~$f$\version{, as in Figure~\ref{fig:cantor}}{}.
	\version{\begin{figure}
		\includegraphics[height=1.5in]{Cantor_function.jpg}
		\includegraphics[height=1.5in]{Primitive_Cantor_function.jpg}
		\caption{Varifold (blue, on the right) with fractal
		``boundary''.  Left: Cantor function. Right: a primitive of
		the Cantor function.} \label{fig:cantor}
	\end{figure}}{}
\end{example}

\begin{definition} [Classes $H(p)$ of summability of mean curvature]
	For~$1 \leq p \leq \infty$, the class~$H(p)$ is defined to consist of
	all $m$-di\-men\-sion\-al varifolds in~$\mathbf R^n$ of locally
	bounded first variation, such that, if $p > 1$, then $\| \updelta V
	\|$ is absolutely continuous with respect to~$\| V \|$, if $1 < p <
	\infty$, their mean curvature is locally $p$-th power $\| V
	\|$~summable, and, if $p = \infty$, it is locally $\| V
	\|$~essentially bounded.
\end{definition}

Observe that $( \int_Z | \mathbf h (V,z) |^p \ud \| V \| \, z )^{1/p}$ is
lower semicontinuous in $V \in H(p)$ if $Z$~is an open subset of~$\mathbf R^n$
and $1 < p < \infty$; a similar statement holds for $p = \infty$.
By Theorem~\ref{thm:integral_compactness}, this entails further compactness
theorems for integral varifolds.

\begin{example} [Critical scaling for $p = m$\version{}{,
see~\protect{\cite[4.2]{MR0307015}}}] \label{example:critical_scaling}

	If $V \in H (m)$, $0 < r < \infty$, and $\phi : \mathbf R^n \to
	\mathbf R^n$ satisfies $\phi (z) = r z$ for $z \in \mathbf R^n$, then
	$\phi_\# V \in H (m)$ and we have $ \int_{\image ( \phi | B )} |
	\mathbf h ( \phi_\# V, z ) |^m \ud \| \phi_\# V \| \, z = \int_B |
	\mathbf h ( V,z ) |^m \ud \| V \| \, z$, and, in case $m=1$, also $\|
	\updelta ( \phi_\# V ) \| ( \image ( \phi | B ) ) = \| \updelta V \| (
	B)$, whenever $B$~is a Borel subset of~$\mathbf R^n$.
\end{example}

\section{Regularity by first variation bounds}

For twice continuously differentiable submanifolds~$M$, the mean curvature
vector $\mathbf h (M,z)$ at~$z \in M$ is defined as trace of the second
fundamental form
\begin{equation*}
	\mathbf b (M,z) : \Tan (M,z) \times \Tan (M,z) \to \Nor (M,z).
\end{equation*}
For varifolds, bounds on the mean curvature (or, more generally, on the first
variation) are defined without reference to a second order structure and
entail more regularity -- similar to the case of weak solutions of the Poisson
equation.  The key additional challenges are non-graphical behaviour and
higher multiplicity:

\begin{example} [Cloud of spheres\version{}{,
see~\protect{\cite[14.1]{MR3528825}}}] \label{example:cloud_of_spheres}

	If $1 \leq p < m < n$ and $Z$ is an open subset of $\mathbf R^n$,
	then, there exists a countable collection~$C$ of $m$-di\-men\-sion\-al
	spheres in~$\mathbf R^n$ such that $V = \sum_{M \in C} \mathbf v_m (C)
	\in H(p)$ and $\spt \| V \| = \Clos Z$.
\end{example}

In contrast, if $p \geq m$, then $\mathscr H^m \restrict \spt \| V \| \leq \|
V \|$, whenever $V$ is an integral varifold in $H(p)$;%
\version{ }{\begin{footnote}%
	{This follows combining \cite[2.7\,(2b), 3.5\,(1c), 8.3]{MR0307015}.}
\end{footnote}}%
in particular, $\spt \| V \|$ has locally finite $m$-di\-men\-sion\-al
measure.

\begin{definition} [Singular]
	An $m$-di\-men\-sion\-al varifold $V$ in~$\mathbf R^n$ is called
	\emph{singular} at~$z$ in~$\spt \| V \|$ if and only if there is no
	neighbourhood of $z$ in which $V$ corresponds to a positive multiple
	of an $m$-di\-men\-sion\-al continuously differentiable submanifold.
\end{definition}

\begin{example} [Zero sets of smooth functions] \label{example:zero_set}
	If $A$~is a closed subset of~$\mathbf R^m$, then there exists a
	nonnegative smooth function $f : \mathbf R^m \to \mathbf R$ with $A =
	\{ x \with f(x)=0 \}$.%
	\version{ }{\begin{footnote}%
		{This assertion may be reduced to the case $m=1$ and $A =
		\mathbf R \cap \{ x \with x \leq 0 \}$
		by~\cite[3.6.1]{MR1014685}.}
	\end{footnote}}%
\end{example}

\begin{example} [Touching\version{}{,
see~\protect{\cite[8.1\,(2)]{MR0307015}}}] \label{example:touching}

	Suppose $A$ is a closed subset of $\mathbf R^m$ without interior
	points such that $\mathscr H^m ( A ) > 0$, $f$ is related to $A$ as in
	Example~\ref{example:zero_set}, $n=m+1$, $\mathbf R^n \simeq \mathbf
	R^m \times \mathbf R$, and $M_i \subset \mathbf R^n$, for $i \in
	\{1,2\}$, correspond to the graph of~$f$, and $\mathbf R^m \times \{ 0
	\}$, respectively.  Then, $V = \mathbf v_m (M_1) + \mathbf v_m (M_2)
	$, see Examples~\ref{example:sum}
	and~\ref{example:varifold_mean_curv}, is an integral varifold in $H (
	\infty )$ which is singular at each point of $M_1 \cap M_2 \simeq A
	\times \{ 0 \}$.
\end{example}

\begin{proposition} [Structure of one-dimensional integral varifolds in
$H(1)$, see~\protect{\cite[12.5]{kol-men.decay}}] \label{thm:1d}

	If $m=1$ and $V \in H(1)$ is an integral varifold, then, near $\| V
	\|$~almost all points, $V$~may be locally represented by finitely many
	Lipschitzian graphs.
\end{proposition}

\begin{proposition} [Regularity of one-dimensional integral varifolds with
vanishing first variation, see\version{ Allard and Almgren
}{~}\protect{\cite[\S\,3]{MR0425741}}] \label{thm:allard_almgren}

	If $m=1$ and $V$ is an integral varifold with $\updelta V = 0$, then
	the set of singular points of $V$ has $\mathscr H^m$ measure zero.
\end{proposition}

It is a key open problem to determine whether the hypothesis $m=1$ is
essential; the difficulty of the case $m \geq 2$ arises from the possible
presence of holes:

\begin{example} [Swiss cheese,
see\version{ Brakke}{~\protect{\cite[6.1]{MR485012}}}
or~\protect{\cite[10.8]{kol-men.decay}}] \label{example:brakke}

	If $2 \leq m < n$, then there exist $V \in H (\infty)$ integral, $T
	\in \mathbf G (n,m)$, and a Borel subset $B$ of $T \cap \spt \| V \|$
	with $\mathscr H^m ( B ) > 0$ such that no $b \in B$ belongs to the
	interior of the orthogonal projection of~$\spt \| V \|$ onto~$T$
	relative to $T$; so, for $\mathscr H^m$ almost all $b \in B$, $\tau
	(b) = T$ and $V$ is singular at~$b$ (cf.\
	Example~\ref{example:sum})\version{; see Figure~\ref{fig:brakke}}{}.
	\version{\begin{figure}
		\includegraphics[height=1.5in]{Catenoids.jpg}
		\includegraphics[height=1.5in]{Plane_with_holes.jpg}
		\caption{Constructing an integral varifold with bounded
		mean curvature whose set of singular points has positive
		weight measure. Left: building blocks (from side). Right:
		result (from top).} \label{fig:brakke}
	\end{figure}}{}
\end{example}

Yet, the preceding varifold needs to have non-singular points in $\spt \| V
\|$.  To treat $\| V \|$~almost all points, however, weaker concepts of
regularity are needed.

\begin{theorem} [\version{Allard's regularity theorem -- qualitative version,
see\version{ Allard }{~}\protect{\cite[8.1\,(1)]{MR0307015}}}{Regularity,
see~\protect{\cite[8.1\,(1)]{MR0307015}}}] \label{thm:allards_regularity}

	If $p > m \geq 2$ and $V \in H(p)$ is integral, then the set of
	singular points of $V$ contains no interior points relative to $\spt
	\| V \|$.
\end{theorem}

\version{\begin{figure}
	\includegraphics[height=1.5in]{Ulrich_Menne.jpg}
	\caption{Ulrich Menne's research concerns geometric measure theory,
	elliptic partial differential equations, and convex analysis.}
\end{figure}}{}

\begin{theorem} [Existence of a second order structure,
see~\protect{\cite[4.8]{MR3023856}}] \label{thm:second_order_rectifiability}

	If $V$ is an integral varifold in $H(1)$, then there exists a
	countable collection~$C$ of \emph{twice} continuously differentiable
	$m$-di\-men\-sion\-al submanifolds such that $C$ covers $\| V
	\|$~almost all of~$\mathbf R^n$.  Moreover, for $M \in C$, there holds
	$\mathbf h ( V,z) = \mathbf h (M,z)$ for $\| V \|$ almost all~$z$.
\end{theorem}

Hence, one may define a concept of second fundamental form of~$V$ such that,
for $\| V \|$~almost all~$z$, its trace is $\mathbf h (V,z)$.  Moreover,
Theorem~\ref{thm:second_order_rectifiability} is the base for studying the
fine behaviour of integral varifolds in~$H(p)$, see \version{\cite{MR3023856}
and \cite{kol-men.decay}}{\cite[\S\,5]{MR3023856} and \cite[\S\,9, \S\,11,
\S\,12]{kol-men.decay}}.


\version{}
{}
\end{document}